\documentclass[12pt,twoside,reqno]{amsart}
\usepackage[colorlinks=true,citecolor=blue]{hyperref}
\usepackage{mathptmx, amsmath, amssymb, amsfonts, amsthm, mathptmx, enumerate, color}
\usepackage{graphicx}
\usepackage{epstopdf}

\theoremstyle{definition}

\numberwithin{equation}{section}

\begin{document}
\setcounter{page}{1}

\vspace*{1.0cm}
\title[ THE LINEAR CONVERGENCE ANALYSIS OF P-GADMM]
{The  global linear convergence rate of  the proximal version of the   generalized alternating direction method of multipliers for separable convex programming}
\author[J. W. Peng, D. X. Liu, X. Q. Zhang, J. C. Yao]{ Jianwen Peng$^{1}$, Dexi Liu$^{1}$, Xueqing Zhang$^{2,*}$, and Jen-Chih Yao$^{3}$}
\maketitle
\vspace*{-0.6cm}

\begin{center}
{\footnotesize {\it

$^1$School of Mathematical Sciences, Chongqing Normal University, Chongqing 400047, P.R.China.\\
$^2$School of Economics and Business Administration, Chongqing University, Chongqing 400044, P.R.China.\\
$^3$Center for General Education, China Medical University,  Taichung 40402, Taiwan.
}}\end{center}

\vskip 4mm {\small\noindent {\bf Abstract.}
To solve the  separable convex optimization problem with linear constraints,  Fang et. al proposed the proximal   version of generalized alternating direction method of multipliers (in short, P-GADMM) and  the doubly  proximal   version of generalized alternating direction method of multipliers (in short, DP-GADMM).  Fang et. al also   proved the worst-case $\mathcal{O}(1/t)$ convergence rate  of both the (P-GADMM) and the (DP-GADMM) measured by the iteration complexity in both   ergodic and   nonergodic senses.  In this paper,   establish the global linear convergence rate   of both the P-GADMM  and the DP-GADMM   for    separable convex optimization problem with the condition  that the subdifferentials of the underlying functions are piecewise linear multifunctions.  We  also illustrate the effectiveness of the   DP-GADMM  by some numerical example about calibrating the correlation matrices.  The results in this paper extend, generalize and improve some known results in the literature.

\noindent {\bf Keywords.}
{Global linear convergence rate; The proximal version of generalized alternating direction method of multipliers; Piecewise linear multifunction; The double  proximal  version of the generalized alternating direction method of multipliers;  The  separable convex optimization problem.
}

\renewcommand{\thefootnote}{}
\footnotetext{ $^*$Corresponding author.
\par
E-mail addresses: jwpeng168@hotmail.com (J. W. Peng), ldx15310888403@163.com (D. X. Liu), zxqcqspb@163.com (X. Q. Zhang), and jen-chih.yao@zjnu.edu.cn (J. C. Yao).
\par
Received January xx, xxxx; Accepted February xx, xxxx. }

\section{Introduction}
In this paper, we consider the following convex minimization model with linear constraints and an objective function which is sum of two functions without coupled variables\\
\begin{equation}
\mathrm{min}~~~~ \{f(x)+g(y)\}~~~~\mathrm{s.t.}~~~~~~Ax+By=b, x\in\mathcal{X},y\in\mathcal{Y},
\end{equation}
where $A\in\mathcal{R}^{\ell\times n}$ and $B\in\mathcal{R}^{\ell\times m}$ are two given matrices, $\mathcal{X}\subseteq\mathcal{R}^{n}$ and $\mathcal{Y}\subseteq\mathcal{R}^{m}$ are two polyhedra, and $f:\mathcal{R}^{n}\rightarrow\mathcal{R}$ and $g:\mathcal{R}^{m}\rightarrow\mathcal{R}$ are proper convex but not necessarily smooth fuctions.\\
\indent   It is well  known that  the classical alternating direction method of multipliers (ADMM) proposed originally in [1,2] is often used to solve the problem (1.1), which is essentially a splitting version of the augmented Lagrangian method in [3, 4]. The separable structure of (1.1) makes it eligible for the application of the ADMM, which reads as\\
\begin{equation}
\begin{split}
\left\{
\begin{aligned}
x^{k+1} &= \mathrm{argmin}\{f(x)-x^{\top}A^{\top}\lambda^{k}+\frac{\beta}{2}\| Ax+By^{k}-b \|^{2}|x\in\mathcal{X}\} , \\
y^{k+1} &= \mathrm{argmin}\{g(y)-y^{\top}B^{\top}\lambda^{k}+\frac{\beta}{2}\| Ax^{k+1}+By-b \|^{2}|y\in\mathcal{Y}\} , \\
\lambda^{k+1}&=\lambda^{k}-\beta(Ax^{k+1}+By^{k+1}-b),\\
\end{aligned}
\right.
\end{split}
\end{equation}
where $\lambda$ is the Lagrange multiplier and $\beta>0$ is a penalty parameter.

 As an efficient and simple acceleration scheme of ADMM,  Eckstein and Bertsekas [5] introduced the following generalized alternating direction method of multipliers (in short, GADMM) whose iterative scheme is \\
\begin{equation}
\begin{split}
\left\{
\begin{aligned}
x^{k+1} =& \mathrm{argmin}\{f(x)-x^{\top}A^{\top}\lambda^{k}+\frac{\beta}{2}\| Ax+By^{k}-b \|^{2}|x\in\mathcal{X}\} , \\
y^{k+1} = &\mathrm{argmin}\{g(y)-y^{\top}B^{\top}\lambda^{k}+\frac{\beta}{2}\| \alpha Ax^{k+1}-(1-\alpha)(By^{k}-b)+By-b \|^{2}|y\in\mathcal{Y}\} , \\
\lambda^{k+1}=&\lambda^{k}-\beta(\alpha Ax^{k+1}-(1-\alpha)(By^{k}-b)+By^{k+1}-b),\\
\end{aligned}
\right.
\end{split}
\end{equation}
where the parameter $\alpha\in(0,2)$ is an acceleration factor. It is easy to see that the original ADMM scheme (1.2) is a special case of the GADMM (1.3) with $\alpha=1$.

In this paper, we consider the  proximal version of generalized alternating direction method of multipliers  (in short, P-GADMM) for (1.1), where only the $x$-subproblems is  proximal regularized,   whose iterative scheme is written as\\
\begin{equation}
\begin{split}
\left\{
\begin{aligned}
x^{k+1}=& \mathrm{argmin}\{f(x)-x^{\top}A^{\top}\lambda^{k}+\frac{\beta}{2}\| Ax+By^{k}-b \|^{2}+\frac{1}{2}\|x-x^{k}\|_{G}^{2}|x\in\chi\} , \\
y^{k+1}=&\mathrm{argmin}\{g(y)-y^{\top}B^{\top}\lambda^{k}+\frac{\beta}{2}\|
\alpha Ax^{k+1}+(1-\alpha)(b-By^{k})+By-b \|^{2}|y\in\mathcal{Y}\}, \\
\lambda^{k+1}=&\lambda^{k}-\beta( \alpha Ax^{k+1}+(1-\alpha)(b-By^{k})+By^{k+1}-b),\\
\end{aligned}
\right.
\end{split}
\end{equation}
where $G$ is a positive definite matrix, $\alpha\in(0,2)$ is an acceleration factor and $\beta>0$ is a penalty parameter.

 It is easy to see that the  proximal version of   alternating direction method of multipliers (in short, P-ADMM)  (i.e.,  (1.6) in [7])  is a special case of the P-GADMM (1.4) with $\alpha=1$.

\indent  When both the $x$- and $y$-subproblems in (1.3) are proximal regularized, the doubly proximal version of the GADMM (in short, DP-GADMM )  for the convex optimization problem (1.1) can be obtained as follows\\
\begin{equation}
\small
\begin{split}
\left\{
\begin{aligned}
x^{k+1}=& \mathrm{argmin}\{f(x)-x^{\top}A^{\top}\lambda^{k}+\frac{\beta}{2}\| Ax+By^{k}-b \|^{2}+\frac{1}{2}\|x-x^{k}\|_{G_{1}}^{2}|x\in\chi\} , \\
y^{k+1}=&\mathrm{argmin}\{g(y)-y^{\top}B^{\top}\lambda^{k}+\frac{\beta}{2}\|
\alpha Ax^{k+1}+(1-\alpha)(b-By^{k})+By-b \|^{2}+\frac{1}{2}\|y-y^{k}\|_{G_{2}}^{2}|y\in\mathcal{Y}\}, \\
\lambda^{k+1}=&\lambda^{k}-\beta( \alpha Ax^{k+1}+(1-\alpha)(b-By^{k})+By^{k+1}-b),\\
\end{aligned}
\right.
\end{split}
\end{equation}
where $G_{i}$ is a positive definite matrix ($i=1, 2$), $\alpha\in(0,2)$ is an acceleration factor and $\beta>0$ is a penalty parameter.

It is noted that the (P-GADMM) and the  (DP-GADMM), respectively,  are called    the linearized version of generalized alternating direction method of multipliers and the doubly  linearized version of generalized alternating direction method of multipliers  by Fang et.al in [6].

 If $\alpha=1$, then DP-GADMM (1.5) becomes the doubly proximal version of the ADMM (in short, DP-ADMM )  for the convex optimization problem (1.1)  as follows\\
\begin{equation}
\begin{split}
\left\{
\begin{aligned}
x^{k+1}=& \mathrm{argmin}\{f(x)-x^{\top}A^{\top}\lambda^{k}+\frac{\beta}{2}\| Ax+By^{k}-b \|^{2}+\frac{1}{2}\|x-x^{k}\|_{G_{1}}^{2}|x\in\chi\} , \\
y^{k+1}=&\mathrm{argmin}\{g(y)-y^{\top}B^{\top}\lambda^{k}+\frac{\beta}{2}\|
  Ax^{k+1}+By-b \|^{2}+\frac{1}{2}\|y-y^{k}\|_{G_{2}}^{2}|y\in\mathcal{Y}\}, \\
\lambda^{k+1}=&\lambda^{k}-\beta(   Ax^{k+1}+By^{k+1}-b),\\
\end{aligned}
\right.
\end{split}
\end{equation}

 It is worthy note that if $G_1=I_{n}$ and  $G_2=I_{m}$, then (1.6) becomes the Proximal ADMM defined by (32) in [8].

The global convergence of ADMM and its variation  was well established under very mild conditions such that the nonemptiness of the solution set of (1.1)  (see [8] and the references therein). However, research on ADMM's iteration complexity   is still in its infancy. Luo and Tseng [9] examined the convergence analysis of matrix splitting algorithms for the affine variational inquality problem on the basis of error bound. Recently,
 Cai and Han [10]    established a worst-case $\mathcal{O}(1/t)$  convergence rate    of  the GADMM. By using the original $\epsilon$-optimal solution measure,
 He and Yuan [7] established the convergence rate $\mathcal{O}(1/n)$ of the P-ADMM in the
worst case   in an ergodic sense, where $n$ denotes the iteration number.  Fang et. al [6] established the worst-case $\mathcal{O}(1/t)$ convergence rate  of both the (P-GADMM) and the (DP-GADMM) measured by the iteration complexity in both   ergodic and   nonergodic senses.
Hong and Luo [11] showed the linear convergence rate of the ADMM for minimizing the sum of finite number of convex separable functions, under a certain error bound condition   and the  sufficiently small dual step size. He, Tao and Yuan [12] established the local linear convergence rate of the ADMM with a backward or forward substitution procedure for solving a separable convex minimization model. When one/some of underlying functions is/are strongly convex, the global linear convergence rate of ADMM was shown in [13,14]. Han, Sun and Zhang [15] estalished the global linear convergence rate of   a semi-proximal ADMM for solving linearly constrained convex composite optimization problems under a mild calmness condition.  Tao and Yuan [16] proved the global linear convergence rate of ADMM for the quadratic program problem.
  Boley [17] and  Han and Yuan [18] proved the local linear convergence rate of GADMM applying to the linear program or  quadratic program.
 Yang and Han [8] established the global linear  convergence rate for both the ADMM and the proximal version of ADMM when applying to (1.1)  that the subdifferentials of the underlying functions are piecewise linear multifunctions. Peng and Zhang [21] established both the local  linear convergence rate    and the global linear convergence rate  of GADMM when applying to (1.1)  that the subdifferentials of the underlying functions are piecewise linear multifunctions.

Inspired by the above works, we investigate the global linear convergence  rate of both the  P-GADMM  and DP-GADMM  when appying to the model (1.1) under the condition   that the  subdifferentials of $f$ and $g$    are piecewise linear multifunctions.

\indent The rest is organized as follows. In section 2, we give some notations and definitions and describe a VI formulation of (1.1),
which are useful for our analysis. We establish the global  linear convergence rate of    the P-GADMM in section 3. In section 4, we show the global linear convergence  rate of DP-GADMM.  In section 5, we give some numerical example to illustrate the effectiveness of the   DP-GADMM. We  give some conclusions of this paper in section 6.\\

\section{Preliminaries}
In this section, we summarize some useful preliminares for further analysis.
In this paper, all vectors are column vectors and for any two vectors $x\in\mathcal{R}^{n}$ and $y\in\mathcal{R}^{m}$, we simply use $u=(x,y)$ to denote their adjunction, i.e., $(x,y)$ denotes $(x^{\top},y^{\top})$. $\|x\|_{p}$ denotes the $p$-norm of vector $x$, where $p=1$ or 2; and for $p=2$, we simply denote it as $\|x\|$. For any positive definite matrix $M$, we denote $\|x\|_{M}=\sqrt{x^{\top}Mx}$ as its $M$-norm. Let $M$ be a positive definite matrix in appropriate dimensionality. For any $u\in \mathcal{U}$, we denote dist$_{M}(u,\mathcal{U}^{*})= $min$\{\|u-u^{*}\|_{M}|u^{*}\in \mathcal{U}^{*}\}$, which represents the distance function from $u$ to $\mathcal{U}^{*}$.\\
\indent For a given matrix $A$, its norm is
\begin{equation*}
\|A\|_{p}:=\sup_{x\neq0}\left\lbrace\frac{\|Ax\|_{p}}{\|x\|_{p}}\right\rbrace.
\end{equation*}
Specially, for a symmetric matrix $A$, $\|A\|_{2}$ denotes its spectral norm.\\
\indent  Let  $\mathcal{K}\subseteq\mathcal{R}^{n}$ be a nonempty, closed, and convex set,    then the projection operator   $P_{\mathcal{K}}(\cdot):\mathcal{R}^{n}\mapsto\mathcal{K}$ be defined by
\begin{equation*}
P_{\mathcal{K}}(u)=\mathrm{argmin}_{v\in {\mathcal{K}}}\{\|v-u\|\},~~\forall u\in\mathcal{R}^{n}.
\end{equation*}
Throught this paper, we will assume that  following condition (i) holds true:\\

\textbf{  Condition (i)}:  The optimal solution set of  {problem (1.1)}  is   nonempty.\\

\textbf{Lemma 2.1.}[8] The following { {results}} are equivalent to each other:\\
(i)  $(x^{*}, y^{*})$ is an optimal solution of (1.1)\\
(ii) there { {exist}} $\xi\in\partial f(x^{*}), \zeta\in\partial g(y^{*})$, and $\lambda^{*}\in\mathcal{R}^{\ell}$ such that
\begin{equation}
\begin{split}
\left\{
\begin{array}{ccc}
&\xi-A^{\top}\lambda^{*}\in-N_{\mathcal{X}}(x^{*}),\\
&\zeta-B^{\top}\lambda^{*}\in-N_{\mathcal{Y}}(y^{*}),\\
&Ax^{*}+By^{*}-b=0.
\end{array}
\right.
\end{split}
\end{equation}
\indent
 where   $N_{\mathcal{S}}(z)$  is { {the}} normal cone of $\mathcal{S}$ at $z$,   $\partial f(x)$ is a subdifferential of $f$ at $x$,  $\partial g(y)$ is a subdifferential of $g$ at $y$.\\
(iii)   $0\in e(u^{*},\gamma)$, or equivalently,  $\mathrm{dist}(0,e(u^{*},\gamma))=0,$\\
where \begin{equation*}       
u^{*}:=\left(                 
  \begin{array}{ccc}   
    x^{*}\\  
    y^{*}\\  
    {\lambda}^{*}\\
      \end{array}
\right)
\in { \mathcal{U}:= } \mathcal{X}\times  \mathcal{Y} \times
    \mathcal{R}^{\ell},
u:=\left(                 
  \begin{array}{ccc}   
    x \\  
    y \\  
    {\lambda} \\
      \end{array}
\right)
\in~ \mathcal{U}, \\
\end{equation*}
  the error function $e(u,\gamma)$    is a set-valued map and defined as follows:\\
 \begin{equation}
e(u,\gamma):=\left(                 
  \begin{array}{ccc}   
    e_{\mathcal{X}}(x,\gamma):=x-P_{\mathcal{X}}[x-\gamma(\partial f(x)-A^{\top}\lambda)]\\  
    e_{\mathcal{Y}}(y,\gamma):=y-P_{\mathcal{Y}}[y-\gamma(\partial g(y)-B^{\top}\lambda)]\\  
    e_{\Lambda}(\lambda,\gamma):=\gamma(Ax+By-b)
  \end{array}
\right).
\end{equation}

If $u^*\in \mathcal{U}$  satisfies   (2.2), then  $u^*$ is a KKT point of problem (1.1).
We use $\mathcal{U}^{*}$ to denote the set of  KKT points  of problem (1.1).

 To simplify our notation in the following analysis, Let\\
\begin{equation}
\begin{split}
v=(y,\lambda); \mathcal{V}=\mathcal{Y}\times R^{\ell}; v^{k}=(y^{k},\lambda^{k}), v^{*}=(y^{*},\lambda^{*}); \\ \Omega^{*}=\{(y^{*},\lambda^{*})|(x^{*},y^{*},\lambda^{*})\in\mathcal{U}^{*}\}.
\end{split}
\end{equation}
\indent Now we define a  matrice $ \Gamma_{\alpha}$ as follows: \\
\begin{equation}       
   \Gamma_{\alpha}=\left(                 
  \begin{array}{ccc}   
  G&0&0\\
0&\frac{\beta}{\alpha}B^{\top}B&\frac{1-\alpha}{\alpha}B^{\top}\\  
0&\frac{1-\alpha}{\alpha}B&\frac{1}{\alpha\beta}I_{\ell}\\  
  \end{array}
\right).              
\end{equation}

Specially,   we can define a matrix $ \Gamma$ as follows:\\

\begin{equation}       
   \Gamma:=\Gamma_{1}=\left(                 
  \begin{array}{ccc}   
  G&0&0\\
0& {\beta}B^{\top}B&0\\  
0& 0 &\frac{1}{\beta}I_{\ell}\\  
  \end{array}
\right).              
\end{equation}

\indent Note the positive definiteness of $\Gamma_{\alpha}$ defined above is guaranteed when $\alpha\in(0,2)$ and $B$ is assumed to be full column rank.

\section{linear convergence rate of the P-GADMM}
In the rest of this paper, we always assume that   the subdifferentials of $f$ and $g$ are piecewise linear multifunctions. for more detaied information on piecewise linear multifunctions, please see [19, 20, 21] and the references therein.

In this section, we will establish the global linear rate of convergence for the P-GADMM to solve the optimization problem (1.1).

\textbf{Lemma 3.1.} (see [6, Thm 1]) Let $\alpha\in(0,2)$ and $B$ is assumed to be full columnrank, then we have
\begin{equation*}
\|u^{k}-u^{*}\|_{\Gamma_{\alpha}}^2-\|u^{k+1}-u^{*}\|_{\Gamma_{\alpha}}^2\geq\|u^{k}-u^{k+1}\|_{\Gamma_{0}}^2.
\end{equation*}
Where
\begin{equation}       
u^{k}:=\left(                 
  \begin{array}{ccc}   
    x^{k}\\  
    y^{k}\\  
    \lambda^{k}\\
  \end{array}
\right),
 u^{*}:=\left(                 
  \begin{array}{ccc}   
    x^{*}\\  
    y^{*}\\  
    \lambda^{*}\\
  \end{array}
\right),
  \Gamma_{0}:=\left(                 
  \begin{array}{ccc}   
    G&0&0\\  
    0&\beta\frac{2-\alpha}{\alpha^{2}}B^{\top}B&0\\  
    0&0&\frac{2-\alpha}{\beta\alpha^{2}}I_{\ell}\\
  \end{array}
\right),          
\end{equation}
 $u^{*}\in \mathcal{U}^{*}$ and $\mathrm{\Gamma_{\alpha}}$ is defined in (2.4).\\

\indent    To prove the linear convergence rate of P-GADMM,
we establish the relationship between $\|u^{k}-u^{k+1}\|_{\Gamma_{0}}^2$ and
  $\mathrm{dist}^2(0,e(u^{k+1},1))$ as follows.\\

\textbf{Lemma 3.2.} Let $u^{k}=(x^{k},y^{k},\lambda^{k})^T$ be the sequence generated by the P-GADMM scheme (1.4). Then, there exists $\sigma>0$ such that\\
\begin{equation}
\|u^{k}-u^{k+1}\|_{\Gamma_{0}}^2\geq \sigma\mathrm{dist}^2(0,e(u^{k+1},1)),
\end{equation}
where $\Gamma_{0}$ is defined in (3.1).\\
$proof.$ By (1.4), there exists $\xi\in\partial f(x^{k+1})$ such that\\
\begin{equation}
x^{k+1}=P_{\mathcal{X}}(x^{k+1}-(\xi-A^{\top}\lambda^{k}+\beta A^{\top}(Ax^{k+1}+By^{k}-b)+G(x^{k+1}-x^{k}))).
\end{equation}
It follows from (3.3) that\\
\begin{equation}
\begin{split}
&\mathrm{dist}(0,e_{\mathcal{X}}(x^{k+1},1))\\
=&\mathrm{dist}(x^{k+1},P_{\mathcal{X}}(x^{k+1}-(\partial f(x^{k+1})-A^{\top}\lambda^{k+1})))\\
\leq&\|P_{\mathcal{X}}(x^{k+1}-(\xi-A^{\top}\lambda^{k}+\beta A^{\top}(Ax^{k+1}+By^{k}-b)\\
&+G(x^{k+1}-x^{k})))-P_{\mathcal{X}}(x^{k+1}-(\xi-A^{\top}\lambda^{k+1}))\|\\
\leq&\|A^{\top}(\lambda^{k}-\lambda^{k+1})-\beta A^{\top}(Ax^{k+1}+By^{k}-b)-G(x^{k+1}-x^{k})\|.
\end{split}
\end{equation}
Thanks to $\lambda^{k+1}=\lambda^{k}-\beta( \alpha Ax^{k+1}+(1-\alpha)(b-By^{k})+By^{k+1}-b)$, then we have\\
\begin{equation}
\begin{split}
\beta(Ax^{k+1}+By^{k}-b)=\frac{1}{\alpha}(\lambda^{k}-\lambda^{k+1})+\frac{\beta}{\alpha}(By^{k}-By^{k+1}).
\end{split}
\end{equation}
Thus, we get\\
\begin{equation}
\begin{split}
&\mathrm{dist}(0,e_{\mathcal{X}}(x^{k+1},1))\\
\leq&\|A^{\top}(\lambda^{k}-\lambda^{k+1})-\beta A^{\top}(Ax^{k+1}+By^{k}-b)-G(x^{k+1}-x^{k})\|\\
=&\|(1-\frac{1}{\alpha})A^{\top}(\lambda^{k}-\lambda^{k+1})+\frac{\beta}{\alpha}A^{\top}(By^{k+1}-By^{k})-G(x^{k+1}-x^{k})\|\\
\leq&|1-\frac{1}{\alpha}|\|A^{\top}\|\|\lambda^{k+1}-\lambda^{k}\|+\frac{\beta}{\alpha}\|A^{\top}\|\|By^{k+1}-By^{k}\|+\|G\|\|x^{k+1}-x^{k}\|.
\end{split}
\end{equation}
From (3.6), we have\\
\begin{equation}
\begin{split}
&\mathrm{dist}^{2}(0,e_{\mathcal{X}}(x^{k+1},1))\\
\leq&3(1-\frac{1}{\alpha})^{2}\|A^{\top}\|^{2}\|\lambda^{k+1}-\lambda^{k}\|^{2}\\
&+3\frac{\beta^{2}}{\alpha^{2}}\|A^{\top}\|^{2}\|By^{k+1}-By^{k}\|^{2}+3\|G\|^{2}\|x^{k+1}-x^{k}\|^{2}.
\end{split}
\end{equation}
By (1.4), there exists $\zeta\in\partial g(y^{k+1})$ such that\\
\begin{equation}
y^{k+1}=P_{\mathcal{Y}}(y^{k+1}-(\zeta-B^{\top}\lambda^{k}+\beta B^{\top}(\alpha Ax^{k+1}+(1-\alpha)(b-By^{k})+By^{k+1}-b))).
\end{equation}
Therefore,
\begin{equation}
\begin{split}
&\mathrm{dist}(0,e_{\mathcal{Y}}(y^{k+1},1))\\
=&\mathrm{dist}(y^{k+1},P_{\mathcal{Y}}(y^{k+1}-(\partial g(y^{k+1})-B^{\top}\lambda^{k+1})))\\
\leq&\|P_{\mathcal{Y}}(y^{k+1}-(\zeta-B^{\top}\lambda^{k}+\beta B^{\top}(\alpha Ax^{k+1}+(1-\alpha)(b-By^{k})+By^{k+1}-b)))\\
&-P_{\mathcal{Y}}(y^{k+1}-(\zeta-B^{\top}\lambda^{k+1}))\|\\
\leq&\|B^{\top}(\lambda^{k}-\lambda^{k+1})-\beta B^{\top}(\alpha Ax^{k+1}+(1-\alpha)(b-By^{k})+By^{k+1}-b)\|\\
=&0,
\end{split}
\end{equation}
where the second equality in (3.9) holds due to the update formula of $\lambda^{k}$ in (1.4).

It follows from $e_{\Lambda}(\lambda,1)=Ax+By-b$ that\\
\begin{equation}
\mathrm{dist}(0,e_{\Lambda}({\lambda}^{k+1},1))\leq\|e_{\Lambda}({\lambda}^{k+1},1)\|=\|Ax^{k+1}+By^{k+1}-b\|
\end{equation}
and\\
\begin{equation}
\begin{split}
&\|Ax^{k+1}+By^{k+1}-b\|\\
=&\|\frac{1}{\alpha\beta}(\lambda^{k}-\lambda^{k+1})-\frac{1}{\alpha}(By^{k+1}-By^{k})+(By^{k+1}-By^{k})\|\\
=&\|\frac{1}{\alpha\beta}(\lambda^{k}-\lambda^{k+1})-\frac{1-\alpha}{\alpha}(By^{k+1}-By^{k})\|\\
\leq&\frac{1}{\alpha\beta}\|\lambda^{k+1}-\lambda^{k}\|+\frac{|1-\alpha|}{\alpha}\|By^{k+1}-By^{k}\|.
\end{split}
\end{equation}
By   (3.10) and (3.11),   we can obtain that\\
\begin{equation}
\mathrm{dist}(0,e_{\Lambda}(\lambda^{k+1},1))\leq\frac{1}{\alpha\beta}\|\lambda^{k+1}-\lambda^{k}\|+\frac{|1-\alpha|}{\alpha}\|By^{k+1}-By^{k}\|.
\end{equation}
So, we get\\
\begin{equation}
\mathrm{dist}^{2}(0,e_{\Lambda}(\lambda^{k+1},1))\leq2\frac{1}{\alpha^{2}\beta^{2}}\|\lambda^{k+1}-\lambda^{k}\|^{2}+2\frac{(1-\alpha)^{2}}{\alpha^{2}}\|By^{k+1}-By^{k}\|^{2}.
\end{equation}
It follows from  (3.7), (3.9) and (3.13) that\\
 \begin{equation}
 \begin{split}
&\mathrm{dist}^{2}(0,e(u^{k+1},1))\\
=&\mathrm{dist}^{2}(0,e_{\mathcal{X}}(x^{k+1},1))+\mathrm{dist}^{2}(0,e_{\mathcal{Y}}(y^{k+1},1))
+\mathrm{dist}^{2}(0,e_{\Lambda}(\lambda^{k+1},1))\\
\leq&(|1-\frac{1}{\alpha}|\|A^{\top}\|\|\lambda^{k+1}-\lambda^{k}\|+\frac{\beta}{\alpha}\|A^{\top}\|\|By^{k+1}-By^{k}\|+\|G\|\|x^{k+1}-x^{k}\|)^{2}\\
&+(\frac{1}{\alpha\beta}\|\lambda^{k+1}-\lambda^{k}\|+\frac{|1-\alpha|}{\alpha}\|By^{k+1}-By^{k}\|)^{2}\\
\leq&3\|G\|^{2}\|x^{k+1}-x^{k}\|^{2}+[\frac{\alpha}{\beta(2-\alpha)}\cdot(\frac{3\beta^{2}}{\alpha}\|A^{\top}\|^{2}+\frac{2(1-\alpha)^{2}}{\alpha})]\cdot\frac{\beta(2-\alpha)}{\alpha^2}\|By^{k+1}-By^{k}\|^{2}\\
&+[\frac{\beta\alpha}{2-\alpha}\cdot(\frac{3}{\alpha}(1-\alpha)^{2}\|A^{\top}\|^{2}+\frac{2}{\alpha\beta^{2}})]\cdot\frac{2-\alpha}{\beta\alpha^2}\|\lambda^{k+1}-\lambda^{k}\|^{2}\\
\leq&\varsigma_{1}(\|x^{k+1}-x^{k}\|^{2}+\frac{\beta(2-\alpha)}{\alpha^2}\|By^{k+1}-By^{k}\|^{2}+\frac{2-\alpha}{\beta\alpha^2}\|\lambda^{k+1}-\lambda^{k}\|^{2})
\end{split}
\end{equation}
where \\
\begin{equation*}
\begin{split}
\varsigma_{1}=&\mathrm{max}\{3\|G\|^{2},\frac{\alpha}{\beta(2-\alpha)}\cdot(\frac{3\beta^{2}}{\alpha}\|A^{\top}\|^{2}+\frac{2(1-\alpha)^{2}}{\alpha}),\frac{\beta\alpha}{2-\alpha}\cdot(\frac{3}{\alpha}(1-\alpha)^{2}\|A^{\top}\|^{2}+\frac{2}{\alpha\beta^{2}}) \}
\end{split}
\end{equation*}
Clearly, there exists $\varsigma_{2}\geq\varsigma_{1}$ such that $\varsigma_{2}G-\varsigma_{1}I$ is positive semidefinite. Then we have \\
 \begin{equation}
 \begin{split}
&\mathrm{dist}^{2}(0,e(u^{k+1},1))\\
\leq&\varsigma_{2}(\|x^{k+1}-x^{k}\|_{G}^{2}+\frac{\beta(2-\alpha)}{\alpha^2}\|By^{k+1}-By^{k}\|^{2}+\frac{2-\alpha}{\beta\alpha^2}\|\lambda^{k+1}-\lambda^{k}\|^{2})\\
=&\varsigma_{2}\|u^{k}-u^{k+1}\|_{\Gamma_{0}}^{2},
\end{split}
\end{equation}
where $\Gamma_{0}$ is definited in (3.1). By $\sigma=\frac{1}{\varsigma_{2}}$, the proof of the Lemma 3.2 is complete.

According to Lemma 3.1, we find that the sequence $u^{k}$ generated by the P-GADMM scheme (1.4) converges to a solution point $u^{*}\in\mathcal{U}^{*}$. Based on this, the sequence $u^{k}$ generated by the P-GADMM scheme (1.4) is bounded. In other words, there exists a bounded set $U$ such that $u^{k} \in U$. Considering that $e(u,\gamma)$ is piecewise linear multifunction of $u$ under the
condition that both $\partial f$ and $\partial g$ are piecewise linear, Based on the result proposed by \text{Yang} and \text{Han} in [8], the following lemma is given, which plays an important role in the proof of our main results. \\

\textbf{Lemma 3.3.} Let $U$ be a bounded set in $\mathcal{R}^{n}$. then there exists $\psi>0$ such that\\
\begin{equation*}
\text{dist}(u,\mathcal{U}^{*})\leq \psi \text{dist}(0,e(u,1)),\forall u\in U.
\end{equation*}

\indent  We establish the global linear convergence rate of P-GADMM as follows:\\

\textbf{Theorem 3.1.} Let $\{u^{k}\}$ be the sequence generated by the P-GADMM scheme (1.4). Then there exists $0<\tau<1$ such that\\
\begin{equation}
\mathrm{dist}_{\Gamma_{\alpha}}^2(u^{k+1},\mathcal{U}^{*})\leq\tau\mathrm{dist}_{\Gamma_{\alpha}}^2(u^{k},\mathcal{U}^{*}),
\end{equation}
where  $\tau=\frac{1}{1+\frac{\sigma}{\psi^2\lambda_{max}(\Gamma_{\alpha})}}$.\\
$proof$.
According to Lemma 3.2 and Lemma 3.3, we can obtain that the following inequality holds
\begin{equation}
\|u^{k}-u^{k+1}\|_{\Gamma_0}^2\geq\sigma\mathrm{dist}^2(0,e(u^{k+1},1))\geq\frac{\sigma}{\psi^2}\text{dist}^2(u^{k+1},\mathcal{U}^{*})\geq\frac{\sigma}{\psi^2\lambda_{max}(\Gamma_{\alpha})}\text{dist}_{\Gamma_{\alpha}}^2(u^{k+1},\mathcal{U}^{*}),
\end{equation}
where $\lambda_{max}(\Gamma_{\alpha})$ is the largest eigenvalue of $\Gamma_{\alpha}$.

By selecting $u^{*} \in\mathcal{U}^{*}$, since\\
\begin{equation*}
\mathrm{dist}_{\Gamma_{\alpha}}(u^{k},\mathcal{U}^{*})=\text{min}\{\|u^{k}-u\|_{\Gamma_{\alpha}}|u\in\mathcal{U}^{*}\},
\end{equation*}
by Lemma 3.1, Lemma 3.2 and (3.17), we have
\begin{equation}
\begin{split}
&(1+\frac{\sigma}{\psi^2\lambda_{max}(\Gamma_{\alpha})})\mathrm{dist}_{\Gamma_{\alpha}}^2(u^{k+1},\mathcal{U}^{*})\\
\leq&\|u^{k+1}-u^{*}\|_{\Gamma_{\alpha}}^2+\frac{\sigma}{\psi^2\lambda_{max}(\Gamma_{\alpha})}\mathrm{dist}_{\Gamma_{\alpha}}^2(u^{k+1},\mathcal{U}^{*})\\
\leq&\|u^{k+1}-u^{*}\|_{\Gamma_{\alpha}}^2+\|u^{k+1}-u^{k}\|_{\Gamma_{0}}^2\\
\leq&\|u^{k}-u^{*}\|_{\Gamma_{\alpha}}^2.
\end{split}
\end{equation}
Hence,  for any $u^{*}\in \mathcal{U}^{*}$,  we have
\begin{equation*}
(1+\frac{\sigma}{\psi^2\lambda_{max}(\Gamma_{\alpha})})\mathrm{dist}_{\Gamma_{\alpha}}^2(u^{k+1},\mathcal{U}^{*})
\leq\mathrm{dist}_{\Gamma_{\alpha}}^2(u^{k},\mathcal{U}^{*}).
\end{equation*}
Let $\tau=\frac{1}{1+\frac{\sigma}{\psi^2\lambda_{max}(\Gamma_{\alpha})}}$, then the assertion (3.16) is proved.

Let $\alpha=1$, then by Theorem 3.1, we get the following global linear convergence rate of the P-ADMM scheme (1.6) in [7]:\\

\textbf{Corollary  3.1.} Let $\{u^{k}\}$ be the sequence generated by the P-ADMM scheme (1.6) in [7]. Then there exists $0<\varphi<1$ such that\\
\begin{equation}
\mathrm{dist}_{\Gamma}^2(u^{k+1},\mathcal{U}^{*})\leq\varphi\mathrm{dist}_{\Gamma}^2(u^{k},\mathcal{U}^{*}),
\end{equation}
where $\Gamma$ is defined in (2.5).\\
\textbf{Remark 3.1.} Corollary  3.1 improves Theorem 4.1 in [7] since  the sublinear convergence rate of P-ADMM hsa been changed as   the global linear convergence rate of P-ADMM.\\
\textbf{Remark 3.2.}  Theorem 3.1 improves Theorems 2 and 6 in [6] since the worst-case $\mathcal{O}(1/t)$ convergence rate measured by the iteration complexity in both   ergodic and   nonergodic senses of the P-GADMM  has been changed as   the global linear convergence rate of the P-GADMM.

\section{The linear convergence rate of the DP-GADMM}
In this section, we show the global linear convergence  of DP-GADMM defined by (1.5).

\indent  We define one auxiliary sequence  for the convenience of analysis. More specifically, for the sequence $\{u^{k}\}$ generated by the DP-GADMM (1.5), let\\
  \begin{equation}
\bar{u}^k=\left (                 
  \begin{array}{ccc}   
    \bar{x}^{k}\\  
    \bar{y}^{k}\\  
\bar{\lambda}^{k}\\
  \end{array}\right)
=\left (                 
  \begin{array}{ccc}   
    x^{k+1}\\  
    y^{k+1}\\  
    \lambda^{k}-\beta(Ax^{k+1}+By^{k}-b)
  \end{array}
\right)
\end{equation}
And for notational simplicity, we also define some matrices in the following analysis. Let

\begin{equation}      
 H_{\alpha}=\left(                 
  \begin{array}{ccc}   
    G_{1} &0& 0\\  
    0 & \frac{\beta}{\alpha}B^{\top}B +G_{2}& \frac{1-\alpha}{\alpha}B^{\top}\\  
    0&\frac{1-\alpha}{\alpha}B&\frac{1}{\alpha\beta}I_{\ell}\\
  \end{array}
\right),\\      
 \end{equation}

If $\alpha=1$, then $H_{\alpha}$ reduces to the following matrix\\
\begin{equation}       
\begin{split}
H=\left(                 
  \begin{array}{ccc}   
    G_{1} &0& 0\\  
    0 & \beta B^{\top}B+G_{2}& 0\\  
    0&0&\frac{1}{\beta}I_{\ell}\\
  \end{array}
\right).\\               
\end{split}
\end{equation}
It is easy to see that the matrices  $H$ and   $H_{\alpha}$  are both  positive definite when $\alpha\in(0,2)$.

Let $\{u^{k}\}$ be the sequence generated by the DP-GADMM scheme (1.5) and consider the definition of $\{\bar{u}^k\}$ in (4.1). Fang et. al examine the boundedness of the sequence $\{u^{k}\}$ generated by the DP-GADMM (1.5) in [6, Thm 9]. The following lemma can be easily  proved by combining Lemma 13 and Theorem 9 in [6].\\

\textbf{Lemma 4.1.} Let $\{u^{k}\}$ be the sequence generated by the DP-GADMM scheme (1.5) with $\alpha\in(0,2)$. Then, it holds that
\begin{equation}
\|u^{k}-u^{*}\|_{H_{\alpha}}^2-\|u^{k+1}-u^{*}\|_{H_{\alpha}}^2\geq\|u^{k}-\bar{u}^k\|_{H_{0}}^2,
\end{equation}
where $u^{k}$ and $u^*$ is defined in (3.1), $\bar{u}^k$ is defined in (4.1), $\mathrm{H_{\alpha}}$ is defined  in (4.2), and $H_{0}$ is defined as follows
 \begin{equation}     
   H_{0}:=\left(                 
  \begin{array}{ccc}   
    G_{1}&0&0\\  
    0&G_{2}&0\\  
    0&0&\frac{2-\alpha}{\beta}I_{\ell}\\
  \end{array}
\right)            
\end{equation}

The following result  plays a vital role in establishing the linear convergence rate of the DP-GADMM (1.5).\\

\textbf{Lemma 4.2.} Let $\{u^{k}\}$ be the sequence generated by the DP-GADMM scheme (1.5) and $\alpha\in(0,2)$. Then there exists a parameter $\delta>0$ so that the following inequality holds
\begin{equation}
\|u^{k}-u^{*}\|_{H_{\alpha}}^2-\|u^{k+1}-u^{*}\|_{H_{\alpha}}^2\geq\delta\|u^{k}-u^{k+1}\|_{H_{0}}^2,
\end{equation}
where $\delta=\frac{1}{max\{1+\frac{2\beta(2-\alpha)\|B\|^2}{\lambda_{min}(G_2)}, 2\alpha^2 \}}$.\\
$proof.$ According to the definition of $\lambda^{k+1}$ in (1.5) and   (4.1), we have
\begin{equation}
\lambda^{k}-\lambda^{k+1}=\alpha (\lambda^{k}-\bar{\lambda}^{k})-\beta B(y^k-\bar{y}^{k}),
\end{equation}
By (4.7) and   the triangle inequality, we have
\begin{equation}
\begin{split}
\|\lambda^k-\lambda^{k+1}\|^{2}&=\|\alpha (\lambda^{k}-\bar{\lambda}^{k})-\beta B(y^k-\bar{y}^{k})\|^2\leq(\|\alpha (\lambda^{k}-\bar{\lambda}^{k})\|+\|\beta B(y^k-\bar{y}^{k})\|)^2\\
&\leq 2\alpha^2\|\lambda^k-\bar{\lambda}^k\|^2+2\beta^2\|B(y^k-y^{k+1})\|^2
\end{split}
\end{equation}
It follows from (4.8) and  (4.1) that
\begin{equation}
\begin{split}
&\|u^k-u^{k+1}\|_{H_{0}}^2=\|x^k-x^{k+1}\|_{G_1}^2+\|y^k-y^{k+1}\|_{G_2}^2+\frac{2-\alpha}{\beta}\|\lambda^k-\lambda^{k+1}\|^2\\
&\leq\|x^k-\bar{x}^{k}\|_{G_1}^2+\|y^k-\bar{y}^{k}\|_{G_2}^2+(2\beta(2-\alpha)\|B\|^2)\cdot\|y^k-\bar{y}^{k}\|^2+(\frac{2\alpha^2(2-\alpha)}{\beta})\cdot\|\lambda^k-\bar{\lambda}^{k}\|^2\\
&\leq\|x^k-\bar{x}^{k}\|_{G_1}^2+(1+\frac{2\beta(2-\alpha)\|B\|^2}{\lambda_{min}(G_2)})\cdot\|y^k-\bar{y}^{k}\|_{G_2}^2+(\frac{2\alpha^2(2-\alpha)}{\beta})\cdot\|\lambda^k-\bar{\lambda}^{k}\|^2\\
&=\|x^k-\bar{x}^{k}\|_{G_1}^2+(1+\frac{2\beta(2-\alpha)\|B\|^2}{\lambda_{min}(G_2)})\cdot\|y^k-\bar{y}^{k}\|_{G_2}^2+\frac{2-\alpha}{\beta}\cdot(2\alpha^2)\cdot\|\lambda^k-\bar{\lambda}^{k}\|^2\\
&\leq max\{1+\frac{2\beta(2-\alpha)\|B\|^2}{\lambda_{min}(G_2)}, 2\alpha^2 \}\cdot\|u^k-\bar{u}^k\|_{H_0}^2,
\end{split}
\end{equation}
where $\lambda_{min}(G_2)$ is the minimum eigenvalue of $G_2$.

Let $\delta=\frac{1}{max\{1+\frac{2\beta(2-\alpha)\|B\|^2}{\lambda_{min}(G_2)}, 2\alpha^2 \}}$, then we have
\begin{equation}
\|u^k-\bar{u}^k\|_{H_0}^2\geq \delta\|u^k-u^{k+1}\|_{H_{0}}^2.
\end{equation}
According (4.4) and (4.10), the conclution of lemma 4.2 holds true.

\indent By Lemma 4.2 and similar auguments   in the proof of the Lemma 3.2, we can easily obtain the following result:

{\bf Lemma 4.3}   Let $u^{k}=(x^{k},y^{k},\lambda^{k})^T$ be the sequence generated by the DP-GADMM scheme (1.5). Then,
there exists a parameter $\sigma'>0$ such that
\begin{equation}
\|u^{k}-u^{k+1}\|_{H_{0}}^2\geq \sigma'\mathrm{dist}^2(0,e(u^{k+1},1)),
\end{equation}
where $H_{0}$ is defined in (4.5).\\

By Lemmas 4.2 and 4.3 and by  similar auguments   in the proof of the   of Theorem 3.1, we can obtain the global linear convergence rate of DP-GADMM (1.5) as follows:\\

\textbf{Theorem 4.1.} Let $\{u^{k}\}$ be the sequence generated by the DP-GADMM scheme (1.5). Then there exists $0<\tau'<1$ such that
\begin{equation}
\mathrm{dist}_{H_{\alpha}}^2(u^{k+1},\mathcal{U}^{*})\leq\tau'\mathrm{dist}_{H_{\alpha}}^2(u^{k},\mathcal{U}^{*}).
\end{equation}

Let $\alpha=1$, then by Theorem 4.1, we obtain the following result:\\

\textbf{Corollary  4.1.} Let $\{u^{k}\}$ be the sequence generated by the DP-ADMM scheme (1.6). Then there exists $0<\tau'<1$ such that\\
\begin{equation}
\mathrm{dist}_{H}^2(u^{k+1},\mathcal{U}^{*})\leq\tau'\mathrm{dist}_{H}^2(u^{k},\mathcal{U}^{*}),
\end{equation}
where $H$ is defined in (4.3).\\

\textbf{Remark  4.1.}  Let  $G_1=I_{n}$ and $G_2=I_{m}$, then by Corollary 4.1, we recover Theorem 17 in [8].\\

\textbf{Remark  4.2.}  Theorems 3.1 and 4.1 extend and generalize Theorem 14 in [8] from ADMM case to both the P-ADMM case and DP-GADMM case.\\

\textbf{Remark 4.3.}  Theorem 3.1 and 4.1  generalize and extend   Theorem 5.1 in [21]  from G-ADMM case to  both the P-ADMM case and DP-GADMM case.

\section{Numerical experiments}
In this section, we report some numerical results for the proposed method (1.5) by calibrating the correlation matrices. All experiments were implemented in MATLAB R2022b on a LAPTOP-HH353THN
with an AMD Ryzen 7 4800U at 1.80 GHz and 16.0 GB memory. We consider to solve the following matrix optimization problem

\begin{equation}
\begin{array}{l}
\mathrm{min}~\frac{1}{2}\|X-C\|^{2}+\frac{1}{2}\|Y-C\|^{2}\\
\mathrm{s.t.}   ~~~~X-Y=\textbf{0},\\
~~~~~~~~X\in S_{+}^{n}, Y\in S_{B},
\end{array}
\end{equation}
where
\begin{equation*}
S_{+}^{n}=\{H\in R^{n\times n}|H^{\top}=H,H\succeq\textbf{0}\},~~ S_{B}=\{H\in R^{n\times n}|H^{\top}=H,H_{L}\preceq H\preceq H_{U}\}.
\end{equation*}
Setting $A=I_n$, $B=-I_n$, we can derive the following subproblem by the algorithm (1.5).
\begin{equation}
\left\{
\begin{aligned}
X^{k+1} =& \mathrm{argmin}\{\frac{1}{2}\|X-C\|^{2}-X^{\top}A^{\top}\lambda^{k}+\frac{\beta}{2}\| AX+BY^{k}-b\|^{2}+\frac{1}{2}\|X-X^{k}\|_{G_{1}}^{2}|X\in\mathcal{X}\} , \\
Y^{k+1} =&\mathrm{argmin}\{\frac{1}{2}\|Y-C\|^{2}-Y^{\top}B^{\top}\lambda^{k}+\frac{\beta}{2}\| \alpha AX^{k+1}+(1-\alpha)(b-BY^{k})+BY-b\|^{2}\\
&+\frac{1}{2}\|Y-Y^{k}\|_{G_{2}}^{2}|Y\in\mathcal{Y}\} , \\
\lambda^{k+1}=&\lambda^{k}-\beta(  \alpha AX^{k+1}+(1-\alpha)(b-BY^{k})+BY^{k+1}-b).
\end{aligned}
\right.
\end{equation}
\begin{figure}[!h]
\centering\includegraphics[width=4.1in]{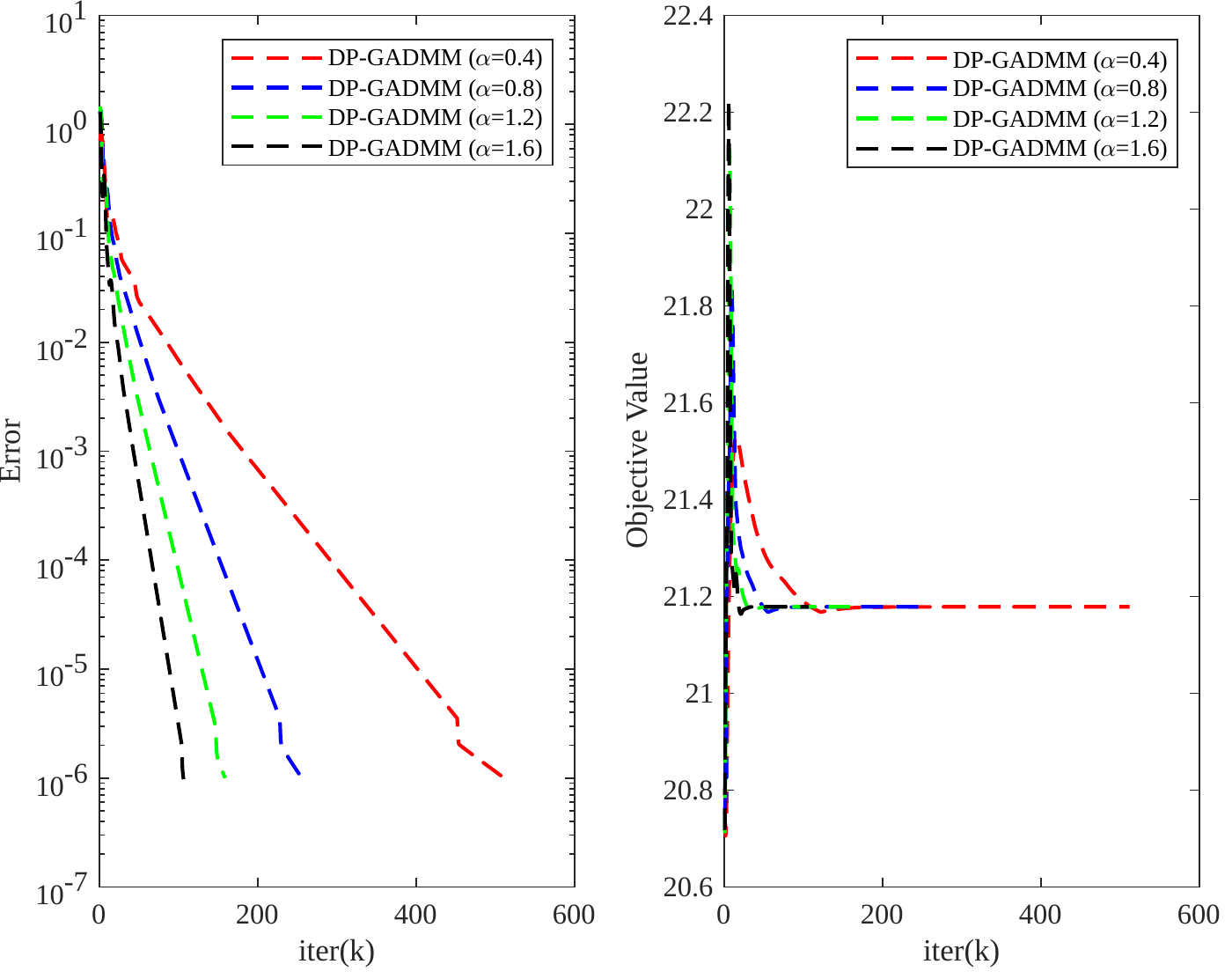}
\caption{The impact of the acceleration factor $\alpha$ ($n=50$)}\label{fig:1}
\end{figure}
In our experiment, we let $C=\mathrm{rand}(n,n)$, $C=(C'+C)-\mathrm{ones}(n,n) + eye(n)$, $HU=\mathrm{ones}(n,n)*0.1$ and $HL=-HU$. And setting $n=50$ and using the DP-GADMM (1.5) to solve the problem (5.1), we can elaborate the convergence rate of the DP-GADMM (1.5) on the basis of different  relaxation factors $\alpha$ (penalty parameter $\beta$), and the stopping criterion is written as
\begin{equation*}
\frac{\mathrm{max}\left\{\|Y^{k}-Y^{k+1}\|,\|\lambda^{k+1}-\lambda^{k}\|\right\}}{\mathrm{max}\left\{\|X^{1}-X^{0}\|,\|Y^{1}-Y^{0}\|,\|\lambda^{1}-\lambda^{0}\|\right\}}<\mathrm{tol}, ~\mathrm{where}~~~ \mathrm{tol}=10^{-6}.
 \end{equation*}
  \begin{figure}[!h]
\centering\includegraphics[width=4.1in]{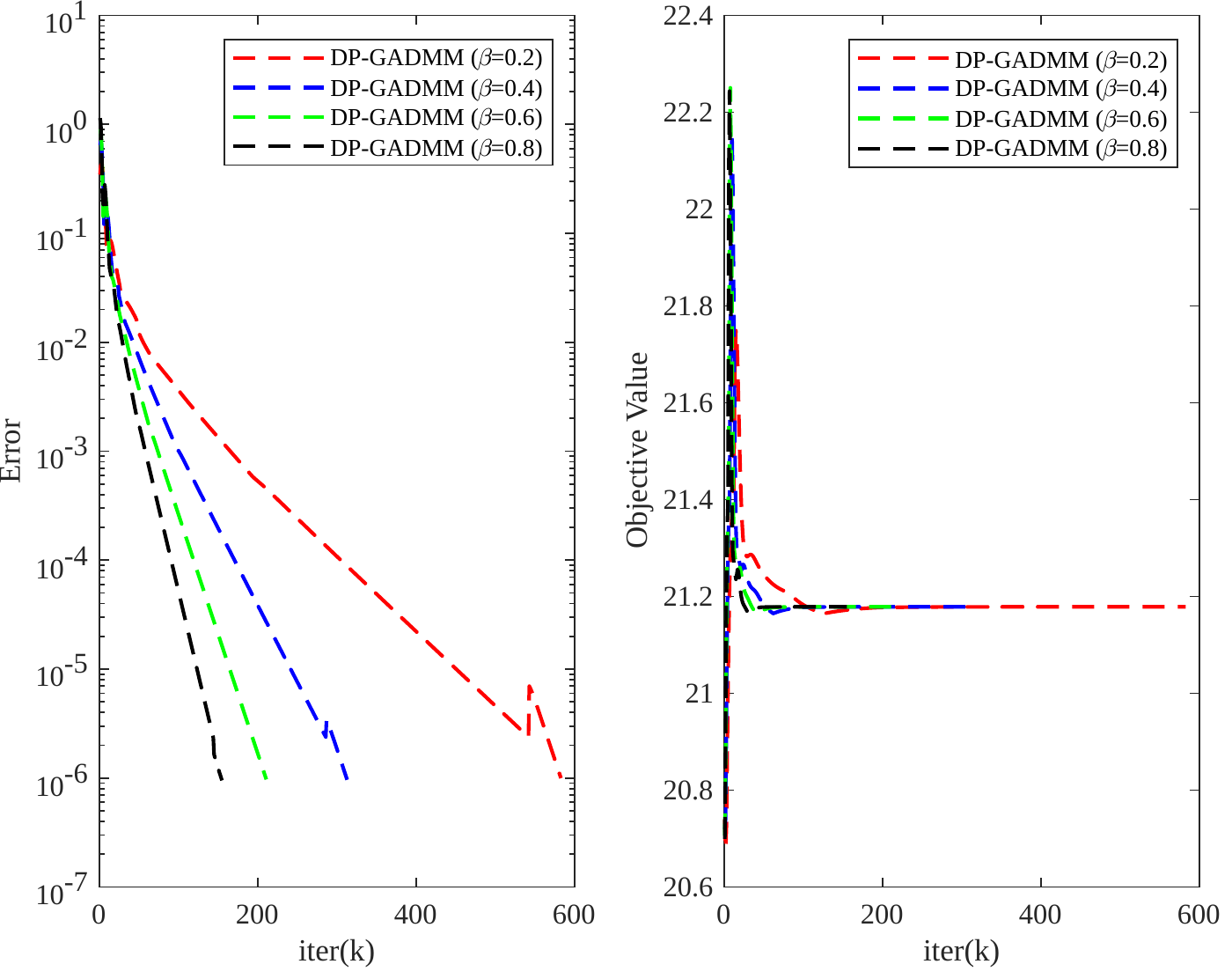}
\caption{The impact of the penalty parameter $\beta$  ($n=50$)}\label{fig:2}
\end{figure}

On  one hand, as shown in Figure 1, we plot the error and number of iterations for the DP-GADMM (1.5). We find that the algorithm (1.5) converges faster as the parameter $\alpha$ increases. In other words, if $\alpha$ is closer   to 2, then  the error is   smaller; if    $\alpha$ is  closer to $0$, the error is greater. And we observe the changes of objective function value with iterations for different values of $\alpha$. We obsearve that  if $\alpha$ is close to 2, then the objective value decreases faster than the cases where $\alpha$ is close to 0.

On the other hand, as shown in Figure 2,   the convergence rate of the  DP-GADMM  accelerates significantly as the Lagrange penalty parameter increases. That is, if the parameter $\beta$ is  the closer to 0, then  the convergence rate is   slower, and  if the parameter $\beta$ is   closer to 1, then  the convergence speedis faster. In addition, we also show a numeric comparison by DP-GADMM, GADMM [5] and ADMM with Gaussian back substitution (ADM-G) [22] by setting $\beta=1$, and $G_{1}=G_{2}=2*I_{n}$ (if $G_{2}=0*I_n$, the  DP-GADMM (1.5) reduces to the P-GADMM (1.4), and the effectiveness of the algorithm (1.4) has been verified numerically by Fang et. al in [6] ). The parameter $\gamma\in(0,1)$ represents step size of the Gaussian back substitution step, we let $\gamma=0.5$ in our experiment. We give these numerical results in Table 1, such as the CPU time, the iteration, Objective value, and Epsilon.
\begin{table}[!h]
\small
\centering  
\begin{tabular}{lccccccccccc}  
\hline
&DP-GADMM $(\alpha=0.5)$ &DP-GADMM $(\alpha=1.5)$&ADM-G &GADMM $(\alpha=0.5)$\\ \hline  
n=50\\         
~~Iterations&412&117&417&379\\
~~CPU times(s)&0.527103&0.183529& 1.011450&1.054792\\        
~~Objective value&21.18&21.18&21.18&21.18\\
~~Epsilon&9.890e-07&9.783e-07&9.972e-07&9.290e-07\\
n=100\\
~~Iterations&571&172&$\gg$1000&$\gg$1000\\
~~CPU times(s)&4.236950 &1.147994&21.687261&19.898802\\        
~~Objective value&47.62&47.62&47.74&47.74\\
~~Epsilon&9.977e-07 &9.721e-07&2.574e-05&2.293e-06\\
n=200\\
~~Iterations&811&290&$\gg$1000&$\gg$1000\\
~~CPU times(s)&24.938693& 8.812811&99.081877&84.167253\\        
~~Objective value&107.55&107.55&107.92&107.96\\
~~Epsilon&9.992e-07&8.180e-07 &2.662e-05&9.999e-07\\\hline
\end{tabular}
\caption{ Numerical comparison between the averaged performance of DP-GADMM and ADM-G} \label{fig:1}
\end{table}
From Table 1, we find that the convergence rate of DP-GADMM highly depend on the parameter $\alpha$. Although setting a large step size ($\gamma$) may yield faster convergence empirically for ADM-G (see [22]), in most cases, DP-GADMM has better convergence than ADM-G. Note that with the increase of $\alpha$, the time spent in  DP-GADMM will be less, and the corresponding number of iterations will also be reduced. However, with the increase of matrix dimension $n$, the number of iterations and CPU time of ADM-G increase significantly. Let $G_1=G_2=0*I_n$, the DP-GADMM (1.5) reduces to the GADMM (1.3). From the experiment, we also find that when setting $\alpha=0.5$, the doubly proximal version of GADMM will converge faster than GADMM, which indicates that the proximal regularization term (e.g., $\frac{1}{2}\|x-x^{k}\|_{G}^{2}$) of algorithm (1.5) has an important influence on the algorithm convergence rate. Therefore, under certain conditions, the convergence rate of DP-GADMM is better than that of ADM-G [22] and GADMM [5].

\section{conclusions}
In this paper, we consider the global  linear convergence rate  of the proximal version of the generalized
alternating direction method of multipliers (P-GADMM)  and the doubly proximal version of the generalized
alternating direction method of multipliers (DP-GADMM)  when optimizing the sum of two
separable convex functions with linear constraints. Based on the error bound, under the condition  that the subdifferentials of the underlying functions are piecewise linear multifunctions, we establish the global linear convergence rate of both the P-GADMM and the DP-GADMM. Finally, we also illustrated the numerical efficiency of the DP-GADMM.

\vskip 6mm
\noindent{\bf Acknowledgments}

\noindent    This research  was partially supported by   the National Natural Science Foundation of China (11991024, 12271071),   the Team Project of Innovation Leading Talent in Chongqing (CQYC20210309536), “Contract System” Project  of Chongqing Talent Plan (cstc2022ycjh-bgzxm0147), the Chongqing University Innovation Research Group Project (CXQT20014) and  the Basic and Advanced Research Project of
Chongqing ( cstc2021jcyj-msxmX0300)


\begin{thebibliography}{99}
\bibitem{[1] }D. Gabay, B. Mercier, A dual algorithm for the solution of nonlinear variational problems via finite element approximations, Comput. Math. Appl., 2(1976), 17-40.
\bibitem{[2] } Glowinski, R., Marrocco, A.: Approximation par $\mathrm{\acute{e}}$l$\mathrm{\acute{e}}$ments finis d'ordre un et r$\mathrm{\acute{e}}$solution parp$\mathrm{\acute{e}}$nalisation-dualit$\mathrm{\acute{e}}$ d'une classe de probl$\mathrm{\acute{e}}$mes non lin$\mathrm{\acute{e}}$aires. R.A.I.R.O., R2, 1975, 41-76.
\bibitem{[3] }Hestenes, M.R.: Multiplier and gradient methods. J. Optim. Theory Appl., 4(1969), 302-320.
\bibitem{[4] }Powell, M.J.D.: A method for nonlinear constraints in minimization problems. In: Fletcher, R. (ed.)
Optimization. Academic Press (1969).
\bibitem{[5] }J. Eckstein and D. P. Bertsekas, On the Douglas-Rachford splitting method and the proximal
point algorithm for maximal monotone operators, Math. Program., 55 (1992), 293-318.
\bibitem{[6] }Ethan X,Fang, B. S. He, Han Liu, Xiaoming Yuan, Generalized alternating direction method of multipliers: new theoretical insights and applications. Math. Prog. Comp., 7(2015), 149-187.
\bibitem{[7] }B. S. He and X. M. Yuan, On the $O(1/n)$ convergence rate of the Douglas-Rachford alternating
direction method, Siam J Numer Anal., 50 (2012), 700-709.
\bibitem{[8] }W. H. Yang and D. Han, linear convergence of the alternating direction method of multipliers
for a class of convex optimization priblems, Siam J Numer Anal., 54 (2016), 625-640.
\bibitem{[9}Z. Q. Luo and P. Tseng, Error bound and convergence analysis of matrix splitting algorithms
for the affine variational inquality problem, Siam J Optim., 2 (1992), 43-54.
\bibitem{[10] } X. J. Cai and D. R. Han, $O(1/t)$ complexity analysis of the generalized
alternating direction method of multipliers, Sci China Math., 62 (2019), 795-808.
\bibitem{[11]} M. Hong and Z. Q. Luo, On the linear convergence of alternating direction method of
multipliers, Math. Program., Ser. A.,  162 (2017), 165-199.

\bibitem{[12]} B.S. He, M. Tao and X.M. Yuan, Convergence rate analysis for the alternating direction
method of multipliers with a substitution procedure for
separable convex programming, Math. Oper. Res.,   42 (2017), 662-691.
\bibitem{[13] }W. Deng and W. T. Yin, On the global and linear convergence of the generalized alternating
direction method of multipliers, J. Sci. Comput. 66 (2016), 889-916.

\bibitem{[14] }
  T.Y.  Lin, S.Q. Ma, and S.Z. Zhang, On the global   linear convergence of the ADMM with
multiblock  variables, SIAM J. Optim., 25 (2015), 1478-1497.

\bibitem{[15]} D.R. Han, D,F. Sun and L.W. Zhang, Linear rate convergence of the alternating
direction method of multipliers for convex
composite programming, Math. Oper. Res.,    43 (2018), 622-637.
\bibitem{[16] } M. Tao and X.M Yuan, On Glowinski's Open Question of Alternating Direction Method of Multipliers, J. Optim. Theory   Appl.,    179(1) (2018),  163-196.
\bibitem{[17] }D. Boley, Local linear convergence of ADMM on quadratic or linear programs, SIAM J. Optim., 23 (2013), 2183-2207.
\bibitem{[18] }D.R. Han and X.M. Yuan, Local linear convergence of the alternating direction method of multipliers
for quadratic programs, Siam J Numer Anal., 51 (2013), 3446-3457.
\bibitem{[19] }X. Y. Zheng and K. F. Ng, Metric subregularity of piecewise linear multifunctions and applications to piecewise linear multiobjective optimization, SIAM J. Optim., 24 (2014),
154-174.

\bibitem{[20] }  S. M. Robinson, Some continuity properties of polyhedral multifunctions, Math. Program.
Stud., 14 (1981), 206-214.

\bibitem{[21] } J. W. Peng and X. Q. Zhang, Linear  convergence rate  of the generalized alternating
direction method of multipliers for a class of convex minimization problems, J Nonlinear Convex A, to appear.
\bibitem{[22] }
B. S. He, M. Tao, X. M. Yuan, Alternating direction method with Gaussian back substitution for
separable convex programming, SIAM J. Optim. 22(2012), 313-340.

\end{thebibliography}
\end{document}